%% file: main.tex
\documentclass[english]{IEEEtran}
\usepackage[T1]{fontenc}
\usepackage[latin9]{inputenc}
\usepackage{color}
\usepackage{verbatim}
\usepackage{amsthm}
\usepackage{amsmath}
\usepackage{amssymb}
\usepackage{graphicx}

\makeatletter
\theoremstyle{plain}
\newtheorem{thm}{\protect\theoremname}
\theoremstyle{definition}
\newtheorem{problem}[thm]{\protect\problemname}
\theoremstyle{plain}
\newtheorem{lem}[thm]{\protect\lemmaname}
\theoremstyle{remark}
\newtheorem{rem}[thm]{\protect\remarkname}
\theoremstyle{plain}
\newtheorem{prop}[thm]{\protect\propositionname}
\theoremstyle{definition}
\newtheorem{defn}[thm]{\protect\definitionname}

\input{preamble.tex}
\usepackage{cite}

\makeatother

\usepackage{babel}
\providecommand{\definitionname}{Definition}
\providecommand{\lemmaname}{Lemma}
\providecommand{\problemname}{Problem}
\providecommand{\propositionname}{Proposition}
\providecommand{\remarkname}{Remark}
\providecommand{\theoremname}{Theorem}

\begin{document}

\title{Data-Driven Allocation of Vaccines for Controlling Epidemic Outbreaks}

\author{Shuo Han, Victor M. Preciado, Cameron Nowzari, and George J. Pappas\allthanks}

\maketitle
\input{abstract.tex}

\section{Introduction\label{sec:Introduction}}

\input{intro.tex}

\section{Preliminaries \& Problem Definition\label{sec:Preliminaries}}

\input{prelim.tex}

\section{Data-Driven Resource Allocation\label{sec:Data-Driven-Resource-Allocation}}

In this section, we develop a mathematical framework to solve the
robust allocation problem described above. Our solution is based on
geometric programming \cite{BKVH07} and its conic extension recently
proposed by Chandrasekaran and Shah in \cite{chandrasekaran2014conic}.
We start our exposition by briefly reviewing some concepts used in
our formulation.

\subsection{Robust Geometric Programming\label{sub:Robust-Geometric-Programming}}

\input{robust_gp.tex}

\subsection{Robust Optimal Resource Allocation\label{sub:Robust-Resource-Allocation}}

\input{robust_allocation.tex}

\subsection{Convex Set of Data-Coherent Networks\label{sub:Uncertainty-set-defined}}

\input{uncertainty_set.tex}

\section{Simulations\label{sec:Simulations}}

\input{simulations.tex}

\section{Conclusions\label{sec:Conclusions}}

\input{conclusions.tex}

\bibliographystyle{ieeetr}
\bibliography{ViralSpread,misc}
%\cleardoublepage{}\newpage{}\pagenumbering{gobble}

\appendix

\subsection{The Saddle Point Theorem\label{sub:The-Saddle-Point-Theorem}}

\input{saddle_pt_thm.tex}

\subsection{Proof of Proposition~\ref{prop:cgp_formulation}\label{sub:Proof-of-Proposition}}

\input{proof_cgp_formulation.tex}
\end{document}

%% file: preamble.tex
\usepackage{url}
\usepackage{mathtools}
\input{affiliation.tex}

\input{symbols.tex}

\newtheorem{assumption}[thm]{Assumption}
\allowdisplaybreaks[3]

%% file: affiliation.tex
\newcommand{\allthanks}{\thanks{The authors are with the Department of Electrical and Systems Engineering, University of Pennsylvania, Philadelphia, PA 19104, USA. {\tt\scriptsize {\{hanshuo,preciado,cnowzari,pappasg\}@seas.upenn.edu}}. This work was supported in part by the NSF under grants CNS-1302222, IIS-1447470, and CNS-1239224, and TerraSwarm, one of six centers of STARnet, a Semiconductor Research Corporation program sponsored by MARCO and DARPA.}}

%% file: symbols.tex
\newcommand{\optmin}{\mathrm{minimize}}
\newcommand{\optmax}{\mathrm{maximize}}
\newcommand{\optst}{\mathrm{subject\; to}}

\newcommand{\optminc}{\mathrm{min.}}

\newcommand{\optstc}{\mathrm{s.t.}}

%% file: abstract.tex
\begin{abstract}
We propose a mathematical framework, based on conic geometric programming,
to control a susceptible-infected-susceptible viral spreading process
taking place in a directed contact network with unknown contact rates.
We assume that we have access to time series data describing the evolution
of the spreading process observed by a collection of sensor nodes
over a finite time interval. We propose a data-driven robust convex
optimization framework to find the optimal allocation of protection
resources (e.g., vaccines and/or antidotes) to eradicate the viral
spread at the fastest possible rate. In contrast to current network
identification heuristics, in which a single network is identified
to explain the observed data, we use available data to define an uncertainty
set containing all networks that are coherent with empirical observations.
Our characterization of this uncertainty set of networks is tractable
in the context of conic geometric programming, recently proposed by
Chandrasekaran and Shah~\cite{chandrasekaran2014conic}, which allows
us to efficiently find the optimal allocation of resources to control
the worst-case spread that can take place in the uncertainty set of
networks. We illustrate our approach in a transportation network from
which we collect partial data about the dynamics of a hypothetical
epidemic outbreak over a finite period of time.%
\begin{comment}
\begin{abstract}
We propose a mathematical framework, based on conic geometric programming,
to control a susceptible-infected-susceptible (SIS) viral spreading
process taking place in a directed contact network of unknown contact
rates. To extract information about these rates, we assume that we
have access to time series describing the evolution of the spreading
process observed from a collection of sensor nodes during a finite
time interval. Using this data series, we propose a data-driven robust
convex optimization framework to find the optimal allocation of protection
resources (e.g., vaccines and/or antidotes) over a set of control
nodes to eradicate the viral spread at the fastest possible rate.
In contrast to current network identification heuristics, in which
a single network is identified to explain the observed data, we use
available data to define an uncertainty set containing all networks
that are coherent with empirical observations. Our characterization
of this uncertainty set of networks is tractable in the context of
\emph{conic geometric programming}, recently proposed by Chandrasekaran
and Shah~\cite{chandrasekaran2014conic}. In this context, we are
able to efficiently find the optimal allocation of resources to control
the worst-case spread that can take place in the uncertainty set of
networks. We illustrate our approach in a transportation network from
which we collect partial data about the dynamics of a hypothetical
epidemic outbreak over a finite period of time.\end{abstract}
\end{comment}
\end{abstract}

%% file: intro.tex
Modeling and analysis of spreading processes in complex networks is
a rich and interdisciplinary research field with a wide range of applications.
Examples include disease propagation in human populations \cite{Bai75,AM91,New02,Vespignani2001PRE,weiss1971asymptotic}
or information spreading in social networks \cite{DM10,EK10,KKT03,LAH07}.
A classical model of disease spreading is the susceptible-infected-susceptible
(SIS) epidemic model \cite{AM91,Bai75}. This model was originally
proposed in the context of `unstructured' populations \cite{weiss1971asymptotic}.
Due to the current availability of accurate datasets describing complex
patterns of network connectivity, the classical SIS model has been
extended to model spreading processes in `networked' populations using
a variety of Markov models \cite{DM10,GMT05,MOK09,MPV02,New02,PJ09,PJ13,PV02,Vespignani2001PRE,WCWF03}.

There is a fast-growing body of literature on containing epidemic
outbreaks given limited control resources. In the context of epidemiology,
these resources can be pharmaceutical (e.g., vaccines and antidotes)
as well as non-pharmaceutical actions (e.g., traffic control and quarantines).
Since these resources are costly, it is of relevance to develop computational
tools to optimize the allocation of resources throughout a population
to control an outbreak. This problem has attracted the attention of
the network science community, resulting in a variety of vaccination
heuristics. For example, Cohen et al. \cite{cohen2003efficient} proposed
a vaccination strategy, called \emph{acquaintance immunization policy,}
and proved it to be much more efficient than random vaccine allocation.
Borgs et al.~\cite{BCGS10} studied theoretical limits in the control
of spreads in undirected network by distributing antidotes. Chung
et at. \cite{chung2009distributing} proposed an immunization strategy
based on PageRank centrality. Similar problems have also been studied
recently in the communication and control community \cite{WRS08,VOK09,DS11,DSV12,VO13,ahn2013global,khanafer2014stability,drakopoulos2014efficient,CP14,mei2014modeling,ramirez2014distributed,NPP14,hayel2014complete}. 

We base our work on \cite{PZEJP13,PZEJP14}, where Preciado et al.
developed a convex optimization framework to find the cost-optimal
distribution of vaccines and antidotes in both directed and undirected
networks. Although current vaccination strategies assume full knowledge
about the network structure and spreading rates, in most practical
applications, this information is only partially known. To elaborate
on this point, let us consider the following setup. Assume that each
node in a network represents subpopulations (e.g., city districts)
connected by edges that are determined by commuting patterns between
districts. In practice, one can use traffic information and geographical
proximity to infer the existence of an edge connecting districts;
however, it is very challenging to use this information to estimate
the contact rates between subpopulations. Inspired by this practical
realization, we consider a networked SIS model taking place in a contact
network with unknown contact rates. To extract information about these
unknown rates, we assume that we have access to time series data describing
the evolution of the spreading process observed by a collection of
sensor nodes over a finite time interval. Such time series data can
be obtained from web services such as Google Flu Trends \cite{GoogleFlu},
or public health agencies such as the Center for Disease Control in
the US \cite{CDCFlu}.

A possible approach to recover the spreading rates is the use of network
identification techniques \cite{materassi2010topological,materassi2009unveiling,gonccalves2008necessary,yuan2011robust,timme2007revealing,nabi2012sieve,6203379,shahrampour2013reconstruction,SP14,FP14}.
However, these techniques are designed to find only one of the many
networks that are coherent with empirical observations \cite{napoletani2008reconstructing}.
Furthermore, as illustrated in \cite{michener1957quantitative,marinazzo2008kernel},
these techniques can lead to unsuccessful network identification.
In contrast to network identification techniques, we propose a data-driven
robust convex optimization framework to find the optimal allocation
of protection resources (e.g., vaccines and/or antidotes) over a set
of control nodes to eradicate the viral spread at the fastest possible
rate. In contrast to current network identification heuristics, in
which a single network is identified to explain the observed data,
we define an uncertainty set containing all networks that are consistent
with the observed data. Our characterization of this uncertainty set
of networks is tractable in the context of \emph{conic geometric programming},
which has recently been proposed by Chandrasekaran and Shah \cite{chandrasekaran2014conic}.
In this context, we are able to efficiently find the optimal allocation
of resources to control the worst-case spread that can take place
in the uncertainty set of networks. We illustrate our approach in
a transportation network from which we collect partial data about
the dynamics of a hypothetical epidemic outbreak over a finite period
of time. We discover that incorporating observations into the uncertainty
set of networks significantly helps reduce the worst-case bound on
the spreading rate. As we increase either the length of time over
which observations are taken or the number of sensor nodes, the bound
on the spreading rate decreases monotonically and converges after
a relatively small number of observations (either in time or in the
number of nodes). Furthermore, even though our allocation algorithm
does not have access to the true underlying contact network, the resulting
allocation performs surprisingly close to the full-knowledge optimal
allocation.

The rest of the paper is organized as follows. In Section~\ref{sec:Preliminaries},
we provide some preliminaries and formulate the problem under consideration.
In Section~\ref{sec:Data-Driven-Resource-Allocation}, we introduce
the conic geometric programming framework and provide the details
about how to cast our problem into this framework. In Section~\ref{sec:Simulations},
we illustrate our approach with numerical simulations using data from
the air transportation network.

%% file: prelim.tex
We begin by introducing the notation and preliminary results needed
in our derivations. In the rest of the paper, we denote by $\mathbb{R}_{+}^{n}$
(respectively, $\mathbb{R}_{++}^{n}$) the set of $n$-dimensional
vectors with nonnegative (respectively, positive) entries. For $d\in\mathbb{N}$,
we define $\left[d\right]$ as the set of integers $\left\{ 1,\ldots,d\right\} $.
We denote vectors using boldface and matrices using capital letters.
We denote by $\mathbf{0}$ the vector of all zeros. Given two vectors
$\mathbf{a}$ and $\mathbf{b}$ of equal dimension, $\mathbf{a}\succeq\mathbf{b}$
indicates component-wise inequality.

\subsection{Graph-Theoretic Nomenclature}

A \emph{weighted}, \emph{directed} graph is defined as the triad $\mathcal{G}\triangleq\left(\mathcal{V},\mathcal{E},\mathcal{W}\right)$,
where $\mathcal{V}\triangleq\left\{ v_{1},\dots,v_{n}\right\} $ is
a set of $n$ nodes, $\mathcal{E}\subseteq\mathcal{V}\times\mathcal{V}$
is a set of ordered pairs of nodes called directed edges, and the
weight function $\mathcal{W}:\mathcal{E}\rightarrow\mathbb{R}_{++}$
associates \textit{positive} real weights to the edges in $\mathcal{E}$.
Throughout the paper, we may use $v_{i}$ and $i$ interchangeably
for all $i\in[n]$. By convention, we say that $\left(v_{j},v_{i}\right)$
is an edge from $v_{j}$ pointing towards $v_{i}$. We define the
in-neighborhood of node $v_{i}$ as $\mathcal{N}_{i}\triangleq\left\{ j\in\left[n\right]:\left(v_{j},v_{i}\right)\in\mathcal{E}\right\} $.
We define the weighted \emph{in-degree} of node $v_{i}$ as $d_{i}\triangleq\sum_{j\in\mathcal{N}_{i}}\mathcal{W}\left((v_{j},v_{i})\right)$.
A directed path from $v_{i_{1}}$ to $v_{i_{l}}$ in $\mathcal{G}$
is an ordered set of vertices $\left(v_{i_{1}},v_{i_{2}},\ldots,v_{i_{l-1}},v_{i_{l}}\right)$
such that $\left(v_{i_{s}},v_{i_{s+1}}\right)\in\mathcal{E}$ for
$s=1,\ldots,l-1$. A directed graph $\mathcal{G}$ is \emph{strongly
connected} if, for every pair of nodes $v_{i},v_{j}\in\mathcal{V}$,
there is a directed path from $v_{i}$ to $v_{j}$. The \emph{adjacency
matrix} of a weighted, directed graph $\mathcal{G}$, denoted by $A_{\mathcal{G}}$,
is an $n\times n$ matrix with entries $a_{ij}=\mathcal{W}\left((v_{j},v_{i})\right)$
if edge $(v_{j},v_{i})\in\mathcal{E}$, and $a_{ij}=0$ otherwise.
In this paper, we only consider graphs with positively weighted edges;
hence, adjacency matrices are always nonnegative. Given an $n\times n$
nonnegative matrix $A$, we can always associate a directed graph
$\mathcal{G}_{A}$ such that $A$ is the adjacency matrix of $\mathcal{G}_{A}$.
Finally, a nonnegative matrix $A$ is \emph{irreducible} if and only
if its associated graph $\mathcal{G}_{A}$ is strongly connected.

Given an $n\times n$ matrix $M$, we denote by $\lambda_{1}\left(M\right),\ldots,\lambda_{n}\left(M\right)$
the eigenvalues of $M$, where we order them according to their magnitudes,
i.e., $\left|\lambda_{1}\right|\geq\left|\lambda_{2}\right|\geq\ldots\geq\left|\lambda_{n}\right|$.
We denote the corresponding eigenvectors by $\mathbf{v}_{1}\left(M\right),\ldots,\mathbf{v}_{n}\left(M\right)$.
We call $\lambda_{1}\left(M\right)$ the spectral radius (or dominant
eigenvalue) of $M$, which we also denote by $\rho\left(M\right)$.

\subsection{\label{sub:Epidemic-Model}SIS Model in Directed Networks}

In our work, we model the spread of a disease using an extension of
the networked discrete-time SIS model proposed by Wang et al. in \cite{WCWF03}.
In contrast to Wang's model, we consider directed networks (instead
of undirected) with non-homogeneous transmission and recovery rates
(instead of homogeneous) as described below. In all SIS models, each
node can be in one out of two possible states: \emph{susceptible}
or \emph{infected}. Over time, nodes switch their states according
to a stochastic process parameterized by (\emph{i}) a set of infection
rates $\left\{ \beta_{ij}\in\left(0,1\right)\right\} _{(v_{j},v_{i})\in\mathcal{E}}$
representing the rates at which an infection can be transmitted through
the edges in the network, and (\emph{ii}) a set of recovery rates
$\left\{ \delta_{i}\in\left(0,1\right)\right\} _{v_{i}\in\mathcal{V}}$
representing the rates at which nodes recover from an infection. We
define $p_{i}(t)$ to be the probability of node $v_{i}$ being infected
at a particular time slot $t\in\mathbb{N}$. In the epidemiological
problem considered herein, it is convenient to associate each node
not to an individual, but a subpopulation living in a particular district%
\footnote{Although in the original networked SIS model \cite{WCWF03}, nodes
represented individuals in a social network, we find the interpretation
of nodes as districts better-suited for epidemiological applications.%
}. In this context, the variable $p_{i}\left(t\right)$ represents
the fraction of the population being infected at time $t$. In the
original model proposed by Wang et al. \cite{WCWF03}, the infection
and recovery rates were assumed to be homogeneous, i.e., $\beta_{ij}=\beta$
and $\delta_{i}=\delta$, and the evolution of $p_{i}\left(t\right)$
was described by a set of difference equations obtained from a mean-field
approximation (see \cite{WCWF03}, eq. (5)--(6)). In our work, we
consider the case of non-homogeneous contact and recovery rates, for
which the set of difference equations can be easily derived to be~\cite{ahn2013global}%
\begin{comment}
\footnote{The set of difference equations herein described is also the result
of a mean-field approximation (see, for example, \cite{MOK09} for
a rigorous analysis of this approximation). Although this approximation
is numerically accurate for many real-world networks, rigorous conditions
for its applicability is a major open problem in this area.%
}
\end{comment}
\begin{multline}
p_{i}\left(t+1\right)=\left(1-p_{i}\left(t\right)\right)\left\{ 1-\prod_{j\in\mathcal{N}_{i}}\left[1-\beta_{ij}p_{j}\left(t\right)\right]\right\} \\
+\left(1-\delta_{i}\right)p_{i}\left(t\right)\label{eq:HeNiSIS dynamics}
\end{multline}
for $i\in\left[n\right]$. This is a system of nonlinear difference
equations for which one can derive sufficient conditions for global
stability as follows. First, notice the following upper bound of \eqref{eq:HeNiSIS dynamics}
\begin{align}
p_{i}\left(t+1\right) & \leq1-\prod_{j\in\mathcal{N}_{i}}\left[1-\beta_{ij}p_{j}\left(t\right)\right]+\left(1-\delta_{i}\right)p_{i}\left(t\right)\nonumber \\
 & \leq\sum_{j\in\mathcal{N}_{i}}\beta_{ij}p_{j}\left(t\right)+\left(1-\delta_{i}\right)p_{i}\left(t\right),\label{eq:Linear Upper Bound}
\end{align}
where the last upper bound is a close approximation of \eqref{eq:HeNiSIS dynamics}
for $p_{i}\left(t\right)\ll1$ and/or $\beta_{i}\ll1$. For convenience,
we define the complementary recovery rate of node $v_{i}$ as $\delta_{i}^{c}\triangleq1-\delta_{i}$,
and the vector $\mathbf{d}^{c}:=\left(\delta_{1}^{c},\ldots,\delta_{n}^{c}\right)^{T}$.
We also define the matrix of infection rates $B_{\mathcal{G}}\triangleq\left[\beta_{ij}\right]$,
where we assume $\beta_{ij}=0$ for all pairs $\left(i,j\right)$
such that $\left(v_{j},v_{i}\right)\notin\mathcal{E}$. Notice that
$B_{\mathcal{G}}$ maintains the same sparsity pattern as $A_{\mathcal{G}}$.
Using the upper bound in \eqref{eq:Linear Upper Bound}, we define
the following linear discrete-time system $\widehat{p}_{i}\left(t+1\right)=\sum_{j\in\mathcal{N}_{i}}\beta_{ij}\widehat{p}_{j}\left(t\right)+\delta_{i}^{c}\widehat{p}_{i}\left(t\right)$,
$i\in\left[n\right]$, which can be written in matrix-vector form
as $\widehat{\mathbf{p}}\left(t+1\right)=M(B_{\mathcal{G}},\mathbf{d}^{c})\widehat{\mathbf{p}}\left(t\right)$,
where $\widehat{\mathbf{p}}\left(t\right)\triangleq\left(\widehat{p}_{1}\left(t\right),\ldots,\widehat{p}_{n}\left(t\right)\right)^{T}$
and the state matrix is given by $M(B_{\mathcal{G}},\mathbf{d}^{c})\triangleq B_{\mathcal{G}}+\mathrm{diag}(\mathbf{d}^{c})$.
Hence, the linear system is asymptotically stable if
\begin{equation}
\rho\left(M(B_{\mathcal{G}},\mathbf{d}^{c})\right)<1,\label{eq:Spectral Condition}
\end{equation}
where the spectral radius $\rho$ of the state matrix $M$ determines
the exponential decay rate of the infection probabilities, i.e., $\left\Vert \widehat{\mathbf{p}}\left(t\right)\right\Vert \leq c\left\Vert \widehat{\mathbf{p}}\left(0\right)\right\Vert \rho^{t}$
for some $c>0$. Since (\ref{eq:Linear Upper Bound}) upper bounds
(\ref{eq:HeNiSIS dynamics}), we have that $\widehat{\mathbf{p}}\left(t\right)\succeq\mathbf{p}\left(t\right)$
for all $t\in\mathbb{N}$ when \textbf{$\widehat{\mathbf{p}}\left(0\right)=\mathbf{p}\left(0\right)$}.
Therefore, the spectral condition in (\ref{eq:Spectral Condition})
is sufficient for global asymptotic stability of the nonlinear model
in (\ref{eq:HeNiSIS dynamics}). Furthermore, the smaller the magnitude
of $\rho\left(M\right)$, the faster the disease dies out.

\subsection{Problem Formulation\label{sub:Problem-Formulation}}

\textbf{}Our main objective is to find the optimal allocation of
control resources to eradicate a disease at the fastest rate possible.
In order to formulate our problem, we first need to describe what
pieces of information are available and what control actions we are
considering. In what follows, we first describe the information available.
In most real epidemiological problems, researchers do not have access
to the spreading rates associated to the links connecting different
districts. Therefore, the exact state matrix $M(B_{\mathcal{G}},\mathbf{d}^{c})$
is usually unknown. In order to extract information about the state
matrix, we consider two different sources of information that are
generally available in epidemiological problems. We classify these
sources as (\emph{i}) \emph{prior information} about the network topology
and parameters of the disease, and (\emph{ii}) \emph{empirical observations}
about the spreading dynamics. In particular, we consider the following
pieces of prior information: 

\begin{enumerate}
\renewcommand{\labelenumi}{P\theenumi.}

\item We assume that the sparsity pattern of the contact matrix $B_{\mathcal{G}}$
is given, although its entries are unknown. This piece of information
may be inferred from geographical proximity, commuting patterns, or
the presence of transportation links connecting subpopulations.

\item We assume that we know the upper and lower bounds on the spreading
rates associated to each edge, i.e., $\beta_{ij}\in\left[\underline{\beta}_{ij},\overline{\beta}_{ij}\right]$,
for all $\left(i,j\right)\in\mathcal{E}$, which may be inferred from
traffic densities and subpopulation sizes.

\item In practice, each district contains a large number of individuals.
Therefore, we can use the average recovery rate in the absence of
vaccination as an estimation of the nodal recovery rate. We denote
this `natural' recovery rate by $\delta_{i}^{0}$, and assume it to
be known.

\end{enumerate}

Apart from these pieces of prior information, we also assume that
we have access to partial observations about the evolution of the
spread over a finite time interval. In particular, we assume that
we observe the dynamics of the disease for $t\in\left[0,T\right]$
from a collection of sensor nodes $\mathcal{V}_{S}\subseteq\mathcal{V}$.
In other words, we have access to the following data set:
\begin{equation}
\mathcal{D}\triangleq\left\{ p_{i}\left(t\right)\colon\mbox{for all }i\in\mathcal{V}_{S},\ t\in\left[T\right]\right\} .\label{eq:Data Series}
\end{equation}
We assume that the data are collected before any control action is
taken; therefore, the evolution of $p_{i}\left(t\right)$ follows
the dynamics in (\ref{eq:HeNiSIS dynamics}) with $\delta_{i}=\delta_{i}^{0}$
(which we assume to be known).

In what follows, we define an uncertainty set that contains all contact
matrices $B_{\mathcal{G}}$ consistent with both empirical observations
and prior knowledge. Based on our prior knowledge described in items
P1--P3 above, we define the following uncertainty set: 
\begin{multline*}
\Delta_{B_{\mathcal{G}}}^{P}\triangleq\{B_{\mathcal{G}}\in\mathbb{R}^{n\times n}\colon\underline{\beta}_{ij}\leq\beta_{ij}\leq\overline{\beta}_{ij},\ \forall(i,j)\in\mathcal{E};\\
\beta_{ij}=0,\ \forall(i,j)\notin\mathcal{E}\}.
\end{multline*}
We also define $\Delta_{B_{\mathcal{G}}}^{D}$ to be the set of contact
matrices that are coherent with the empirical observations $\mathcal{D}$:
\begin{multline*}
\Delta_{B_{\mathcal{G}}}^{D}\triangleq\{B_{\mathcal{G}}\in\mathbb{R}^{n\times n}\colon\left\{ \beta_{ij}\right\} _{\left(i,j\right)\in\mathcal{E}}\mbox{ satisfy \eqref{eq:HeNiSIS dynamics}}\\
\mbox{for }\delta_{i}=\delta_{i}^{0}\mbox{ and }p_{i}\left(t\right)\in\mathcal{D},\ \forall i\in\mathcal{V}_{S},\ t\in\left[T\right]\}.
\end{multline*}
The set contains those contact matrices $B_{\mathcal{G}}$ such that
the transmission rates $\{\beta_{ij}\}$ are consistent with the `natural'
disease dynamics in (\ref{eq:HeNiSIS dynamics}) with $\delta_{i}=\delta_{i}^{0}$.
Notice that~$\Delta_{B_{\mathcal{G}}}^{D}$ is defined as a collection
of polynomial equality constraints on the contact rates $\{\beta_{ij}\}$
given by (\ref{eq:HeNiSIS dynamics}). The uncertainty set that combines
information from both prior knowledge and empirical observations is
defined as
\[
\Delta_{B_{\mathcal{G}}}\triangleq\Delta_{B_{\mathcal{G}}}^{P}\cap\Delta_{B_{\mathcal{G}}}^{D}.
\]

Having introduced the pieces of available information, we now describe
the set of control actions under consideration. In order to eradicate
the disease at the fastest rate possible, we assume that we can use
pharmaceutical resources to tune the recovery rates in a collection
of control nodes, i.e., $\delta_{i}$ for $v_{i}\in\mathcal{V}_{C}\subseteq\mathcal{V}$.
In practice, these resources might be implemented by, for example,
distributing vaccines and/or antidotes throughout the subpopulations
located at those control districts. We assume that distributing vaccines
in a district has an associated cost, which we represent as a node-dependent
vaccine cost function. It is convenient to describe the vaccine cost
function of a district in terms of its complementary recovery rate
$\delta_{i}^{c}$. We denote the vaccine cost function of node $i$
by $g_{i}\left(\delta_{i}^{c}\right)$. This function represents the
cost of tuning the complementary recovery rate of the subpopulation
at node $i\in\mathcal{V}_{C}$ towards the value $\delta_{i}^{c}$.
We assume that we can control the complementary recovery rate~$\delta_{i}^{c}$
within a given feasible interval $\left[\underline{\delta}_{i}^{c},\overline{\delta}_{i}^{c}\right]$,
where $0<\underline{\delta}_{i}^{c}<\overline{\delta}_{i}^{c}=1-\delta_{i}^{0}$.
We assume that the cost of achieving $\overline{\delta}_{i}^{c}$
is zero, since it is equivalent to maintaining the natural recovery
rate. We also assume that the maximum of $g_{i}$ in $\left[\underline{\delta}_{i}^{c},\overline{\delta}_{i}^{c}\right]$
is achieved at $\underline{\delta}_{i}^{c}$. Furthermore, we also
assume that~$g_{i}$ is monotonically decreasing in the range $\left[\underline{\delta}_{i}^{c},\overline{\delta}_{i}^{c}\right]$.
In other words, as we increase the level of investment to protect
a given subpopulation, we also increase the recovery rate of that
subpopulation. 

We are now in a position to state the control problem under consideration:
\begin{problem}
\emph{\label{prob:Main Problem}(Data-driven optimal} \emph{allocation)
Assume we are given the following pieces of information about a viral
spread:}

\emph{(i) prior information about the state matrix (as described in
P1--P3);}

(\emph{ii})\emph{ a finite (and possibly sparse) data series representing
partial evolution of the spread over a set of sensor nodes $\mathcal{V}_{S}\subseteq\mathcal{V}$
during the time interval $t\in\left[T\right]$ (i.e., $\mathcal{D}$
in \eqref{eq:Data Series});}

\emph{(iii) a set of vaccine cost functions $g_{i}$ for all $i\in\mathcal{V}_{C}$,
and a range of feasible recovery rates} $\left[\underline{\delta}_{i}^{c},\overline{\delta}_{i}^{c}\right]$
such that\emph{ $1-\delta_{i}^{0}=\overline{\delta}_{i}^{c}\geq\delta_{i}^{c}\geq\underline{\delta}_{i}^{c}>0$;}

\emph{(iv) a fixed budget $C>0$ to be allocated throughout a set
of control nodes in $\mathcal{V}_{C}\subseteq\mathcal{V}$, so that
$\sum_{i\in\mathcal{V}_{C}}g_{i}(\delta_{i}^{c})\leq C$.}

\emph{Find the cost-constrained allocation of control resources to
eradicate the disease at the fastest possible exponential rate, measured
as $\rho(M(B_{\mathcal{G}},\mathbf{d}^{c}))$, over} \emph{the uncertainty
set $\Delta_{B_{\mathcal{G}}}$ of contact matrices coherent with
prior knowledge and the observations in $\mathcal{D}$.}
\end{problem}
From the perspective of optimization, Problem \ref{prob:Main Problem}
is equivalent to finding the optimal allocation of resources to minimize
the worst-case (i.e., maximum possible) decay rate $\rho(M(B_{\mathcal{G}},\mathbf{d}^{c}))$
for all $B_{\mathcal{G}}\in\Delta_{B_{\mathcal{G}}}$. This can be
cast as a robust optimization problem in the following:
\begin{alignat}{2}
 & \underset{\mathbf{\mathbf{d}^{c}}}{\optmin} & \quad & \sup_{B_{\mathcal{G}}\in\Delta_{B_{\mathcal{G}}}}\rho(M(B_{\mathcal{G}},\mathbf{d}^{c}))\label{eq:prob_robust_alloc}\\
 & \optst & \quad & \sum_{i\in\mathcal{V}_{C}}g_{i}\left(\delta_{i}^{c}\right)\leq C,\nonumber \\
 &  & \quad & \delta_{i}^{c}\leq\delta_{i}^{c}\leq\overline{\delta}_{i}^{c},\quad i\in\mathcal{V}_{C},\nonumber 
\end{alignat}
where the first constraint accounts for our budget limit $C$. In
general, the set $\Delta_{B_{\mathcal{G}}}$ is nonconvex due to the
observation-based uncertainty set $\Delta_{B_{\mathcal{G}}}^{D}$.
In Section~\ref{sub:Uncertainty-set-defined}, we will define a convex
superset $\widehat{\Delta}_{B_{\mathcal{G}}}^{D}\supset\Delta_{B_{\mathcal{G}}}^{D}$,
such that problem~\eqref{eq:prob_robust_alloc} can be relaxed into
a conic geometric program. In our numerical simulations, we verify
that this relaxation provides a good approximation based on real network
data. From here on, we will refer to problem~\eqref{eq:prob_robust_alloc}
as the \emph{robust allocation problem}.

%% file: robust_gp.tex
Geometric programs (GPs) are a type of quasiconvex optimization problem
that can be easily transformed into a convex program and solved in
polynomial time. Let $x_{1},\ldots,x_{n}>0$ denote $n$ decision
variables and define $\mathbf{x}\triangleq\left(x_{1},\ldots,x_{n}\right)\in\mathbb{R}_{++}^{n}$.
In the context of GP, a \emph{monomial $m(\mathbf{x})$} is defined
as a real-valued function of the form $m(\mathbf{x})\triangleq dx_{1}^{a_{1}}x_{2}^{a_{2}}\ldots x_{n}^{a_{n}}$
with $d>0$ and $a_{i}\in\mathbb{R}$. A \emph{posynomial} function
$f(\mathbf{x})$ is defined as a sum of monomials, i.e., $f(\mathbf{x})\triangleq\sum_{k=1}^{K}c_{k}x_{1}^{a_{1k}}x_{2}^{a_{2k}}\ldots x_{n}^{a_{nk}}$,
where $c_{k}>0$ and $a_{ik}\in\mathbb{R}$. It is convenient to write
down a posynomial as the product of a vector of nonnegative coefficients
$\mathbf{c}\triangleq\left(c_{1},\ldots,c_{K}\right)$ and a vector
of monomials $\mathbf{m}\left(\mathbf{x}\right)\triangleq\left(m_{1}\left(\mathbf{x}\right),\ldots,m_{K}\left(\mathbf{x}\right)\right)^{T}$,
such that $f(\mathbf{x})=\mathbf{c}^{T}\mathbf{m}\left(\mathbf{x}\right)$.
Notice that $\left\{ m_{k}\left(\mathbf{x}\right)\right\} _{k=1}^{K}$
is the set of all $K$ monomials involved in our posynomial. Posynomials
are closed under addition, multiplication, and nonnegative scaling.
A posynomial can be divided by a monomial, with the result a posynomial.

A GP is an optimization problem of the form (see \cite{BKVH07} for
a comprehensive treatment):
\begin{alignat}{2}
 & \underset{\mathbf{x}\in\mathbb{R}_{++}^{n}}{\optmin} & \quad & f_{0}(\mathbf{x})\label{eq:General GP}\\
 & \optst & \quad & f_{i}(\mathbf{x})\leq1,\quad i\in\left[m\right],\nonumber \\
 &  &  & h_{j}(\mathbf{x})=1,\quad j\in\left[p\right],\nonumber 
\end{alignat}
where $f_{i}$ are posynomial functions and $h_{j}\left(\mathbf{x}\right)\triangleq d_{j}x_{1}^{b_{1,j}}x_{2}^{b_{2,j}}\ldots x_{n}^{b_{n,j}}$
are monomials. To write~$f_{i}$ in vector-product form, we can define
a vector $\mathbf{c}_{i}$ of positive coefficients such that $f_{i}\left(\mathbf{x}\right)=\mathbf{c}_{i}^{T}\mathbf{m}(\mathbf{x})$,
so that the posynomial constraints in~\eqref{eq:General GP} can
be written as $\mathbf{c}_{i}^{T}\mathbf{m}(\mathbf{x})\leq1$.

A GP is a quasiconvex optimization problem \cite{BV04} that can be
convexified using the logarithmic change of variables $y_{i}=\log x_{i}$
(see \cite{BKVH07} for more details on this transformation). After
this transformation, the GP in (\ref{eq:General GP}) takes the form
\begin{alignat}{2}
 & \underset{\mathbf{y}\in\mathbb{R}^{n}}{\optmin} & \quad & \widetilde{f}_{0}\left(\mathbf{y}\right)\label{eq:Transformed GP}\\
 & \optst & \quad & \widetilde{f}_{i}\left(\mathbf{y}\right)\leq0,\quad i\in\left[m\right],\nonumber \\
 &  &  & \mathbf{b}_{j}^{T}\mathbf{y}+\log d_{j}=0,\quad j\in\left[p\right],\nonumber 
\end{alignat}
where $\widetilde{f}_{i}\left(\mathbf{y}\right)\triangleq\log f_{i}(e^{\mathbf{y}})$
for $i\in\{0,1,\dots,m\}$ and $\mathbf{b}_{j}\triangleq\left(b_{1,j},\dots,b_{n,j}\right)^{T}$
(i.e., the exponents of the monomial $h_{j}$) for $i\in[m]$. As
a result of this transformation, the optimization problem~\eqref{eq:Transformed GP}
is convex and can be efficiently solved in polynomial time (see~\cite[Chapter 4.5]{BV04}
for more details).

In this paper, we shall use conic GP, which is a conic extension of
GP, to solve the following \emph{robust GP with coefficient uncertainties}:
\begin{alignat}{2}
 & \underset{\mathbf{x}\in\mathbb{R}_{++}^{n}}{\optmin} & \quad & f_{0}(\mathbf{x})\label{eq:robust_GP}\\
 & \optst & \quad & \sup_{\mathbf{c}_{i}\in\mathcal{C}_{i}}\mathbf{c}_{i}^{T}\mathbf{m}\left(\mathbf{x}\right)\leq1,\quad i\in\left[m\right],\label{eq:robust_GP_constr}\\
 &  &  & h_{j}(\mathbf{x})=1,\quad j\in\left[p\right],\nonumber 
\end{alignat}
where~$\mathbf{c}_{i}\in\mathbb{R}_{+}^{K}$ is a vector of coefficients
contained in an uncertainty set $\mathcal{C}_{i}\subseteq\mathbb{R}_{+}^{K}$.
The robust GP in (\ref{eq:robust_GP}) extends the formulation of
the standard GP in (\ref{eq:General GP}) to account for uncertainties
in the coefficients of the posynomial functions $f_{i}$ for $i\in\left[m\right]$. 

However, the constraints~\eqref{eq:robust_GP_constr} cannot be handled
naturally by numerical optimization solvers. In what follows, we propose
a methodology to rewrite these constraints in a more numerically favorable
manner when the uncertainty sets~$\mathcal{C}_{i}$ in~\eqref{eq:robust_GP_constr}
can be expressed in terms of an $m_{i}$-dimensional convex cone $\mathcal{K}_{i}\subset\mathbb{R}^{m_{i}}$
as follows: 
\begin{equation}
\mathcal{C}_{i}\triangleq\{\mathbf{c}_{i}\in\mathbb{R}_{+}^{K}\colon F_{i}\mathbf{c}_{i}+\mathbf{g}_{i}\in\mathcal{K}_{i}\}\label{eq:Tractable Set}
\end{equation}
for some fixed $F_{i}\in\mathbb{R}^{m_{i}\times K}$ and $\mathbf{g}_{i}\in\mathbb{R}^{m_{i}}$.
\begin{comment}
For the uncertainty set used in our robust allocation problem, $\mathcal{K}_{i}$
turns out to be the nonnegative orthant. 
\end{comment}
Based on the representation~\eqref{eq:Tractable Set} of $\mathcal{C}_{i}$,
we can use duality theory to derive a more numerically favorable representation
of the constraint in (\ref{eq:robust_GP_constr}) as follows. Assuming
$\mathcal{C}_{i}$ can be represented as (\ref{eq:Tractable Set}),
we have that for each $i\in\left[m\right]$, constraint (\ref{eq:robust_GP_constr})
is equivalent to the optimal value $P_{i}^{*}$ of the following optimization
problem satisfying $P_{i}^{*}\leq1$:
\begin{alignat*}{2}
P_{i}^{*}\triangleq\; & \underset{\mathbf{c}_{i}}{\optmax} & \quad & \mathbf{c}_{i}^{T}\mathbf{m}\\
 & \optst & \quad & F_{i}\mathbf{c}_{i}+\mathbf{g}_{i}\in\mathcal{K}_{i},\\
 &  &  & \mathbf{c}_{i}\succeq\mathbf{0}.
\end{alignat*}
The dual problem of the above is given by
\begin{alignat*}{2}
 & \underset{\boldsymbol{\nu}_{i}}{\optmin} & \quad & \mathbf{g}_{i}^{T}\boldsymbol{\nu}_{i}\\
 & \optst & \quad & F_{i}^{T}\boldsymbol{\nu}_{i}+\mathbf{m}\preceq\mathbf{0},\\
 &  &  & \boldsymbol{\nu}_{i}\in\mathcal{K}^{*},
\end{alignat*}
where $\mathcal{K}^{*}$ is the dual cone of $\mathcal{K}$~\cite{BV04}.
Assume that strong duality holds in this case. Then the optimal value
of the dual problem is also given by $P_{i}^{*}$. Namely, there exists
a dual feasible $\boldsymbol{\nu}_{i}$ such that $\mathbf{g}_{i}^{T}\boldsymbol{\nu}_{i}=P_{i}^{*}$.
Therefore, the constraint in (\ref{eq:robust_GP_constr}) is equivalent
to:
\begin{equation}
\exists\boldsymbol{\nu}_{i}\in\mathcal{K}^{*}\text{ s.t. }F_{i}^{T}\boldsymbol{\nu}_{i}+\mathbf{m}(\mathbf{x})\preceq\mathbf{0},\ \mathbf{g}_{i}^{T}\boldsymbol{\nu}_{i}\leq1\label{eq:robust_GP_constr_dual}
\end{equation}
for each $i\in\left[m\right]$. For the uncertainty set used in our
robust allocation problem, both $\mathcal{K}$ and $\mathcal{K}^{*}$
are the nonnegative orthant. Therefore, we can use the new constraints
in~\eqref{eq:robust_GP_constr_dual} to replace those in \eqref{eq:robust_GP_constr}
and rewrite the robust GP in~\eqref{eq:robust_GP} as
\begin{alignat}{2}
 & \underset{\mathbf{x}\in\mathbb{R}_{++}^{n},\left\{ \boldsymbol{\nu}_{i}\right\} _{i=1}^{m}}{\optmin} & \quad & f_{0}(\mathbf{x})\label{eq:conic_GP}\\
 & \underset{\hphantom{\mathbf{x}\in\mathbb{R}_{++}^{n},\left\{ \boldsymbol{\nu}_{i}\right\} _{i=1}^{m}}}{\optst} & \quad & \boldsymbol{\nu}_{i}\succeq\mathbf{0},\nonumber \\
 &  &  & F_{i}^{T}\boldsymbol{\nu}_{i}+\mathbf{m}(\mathbf{x})\preceq\mathbf{0},\quad\mathbf{g}_{i}^{T}\boldsymbol{\nu}_{i}\leq1,\nonumber \\
 &  &  & h_{j}(\mathbf{x})=1,\nonumber \\
 &  &  & \mbox{for all }i\in\left[m\right],\ j\in\left[p\right].\nonumber 
\end{alignat}
By applying the logarithmic transformation $y_{i}=\log x_{i}$ to
(\ref{eq:conic_GP}) for all $i\in[m]$, we obtain
\begin{alignat}{2}
 & \underset{\mathbf{y}\in\mathbb{R}^{n},\left\{ \boldsymbol{\nu}_{i}\right\} _{i=1}^{m}}{\optmin}\quad &  & \widetilde{f}_{0}(y)\label{eq:conic_GP-2}\\
 & \underset{\hphantom{\mathbf{y}\in\mathbb{R}^{n},\left\{ \boldsymbol{\nu}_{i}\right\} _{i=1}^{m}}}{\optst}\quad &  & \boldsymbol{\nu}_{i}\succeq\mathbf{0},\nonumber \\
 &  &  & F_{i}^{T}\boldsymbol{\nu}_{i}+\widetilde{\mathbf{m}}(\mathbf{y})\preceq\mathbf{0},\quad\mathbf{g}_{i}^{T}\boldsymbol{\nu}_{i}\leq1,\nonumber \\
 &  &  & \mathbf{b}_{j}^{T}\mathbf{y}+\log d_{j}=0,\nonumber \\
 &  &  & \mbox{for all }i\in\left[m\right],\ j\in\left[p\right],\nonumber 
\end{alignat}
where $\widetilde{f}_{0}(\mathbf{y})=f_{0}(\exp\{\mathbf{x}\})$ and
$\widetilde{\mathbf{m}}(\mathbf{y})=\mathbf{m}(\exp\{\mathbf{x}\})$
($\exp\{\mathbf{x}\}$ is component-wise exponential). It can be shown
that both $\widetilde{f}_{0}$ and the entries of $\widetilde{\mathbf{m}}$
are convex in $\mathbf{y}$, since they are nonnegative sums of exponentials
of affine functions in $\mathbf{y}$~\cite{BV04}. In fact, problem~\eqref{eq:conic_GP-2}
is a convex problem and is a particular instance of a \emph{conic
geometric program}~\cite{chandrasekaran2014conic}. Problems in the
form of (\ref{eq:conic_GP-2}) can be solved efficiently using off-the-shelf
software such as CVX~\cite{cvx}.

%% file: robust_allocation.tex
In the following, we show how to formulate the optimization problem
(\ref{eq:prob_robust_alloc}) as a conic GP using the methodology
proposed in Section~\ref{sub:Robust-Geometric-Programming}. In our
derivations, we use the theory of nonnegative matrices and the Perron-Frobenius
lemma~\cite{meyer2000matrix}:
\begin{lem}
[Perron-Frobenius]\label{prop:p-f_thm}Suppose~$M$ is an irreducible
nonnegative matrix. Then, the spectral radius~$\rho(M)$ of $M$
satisfies:

\emph{(}a\emph{)} $\rho\left(M\right)=\lambda_{1}\left(M\right)>0$
is a simple eigenvalue of $M$;

\emph{(}b\emph{)} $M\mathbf{u}=\rho\left(M\right)\mathbf{u}$ for
some $\mathbf{u}\in\mathbb{R}_{++}^{n}$;

\emph{(}c\emph{)} \textup{$\rho(M)=\inf\{\lambda\in\mathbb{R}\colon M\mathbf{u}\preceq\lambda\mathbf{u}\mbox{ for some }\mathbf{u}\succ\mathbf{0}\}$.}\end{lem}
\begin{rem}
Note that the state matrix $B_{\mathcal{G}}+\mathrm{diag}(\mathbf{d}^{c})$
of the linear system \eqref{eq:Linear Upper Bound} is irreducible
if the graph $\mathcal{G}$ is strongly connected. In what follows,
we shall assume that the contact network $\mathcal{G}$ is strongly
connected. This assumption is reasonable in the context of epidemic
control, since the transportation network connecting different districts
or subpopulations is strongly connected in most cases. Notice also
that, as a consequence of this assumption, all the matrices in the
uncertainty set $\Delta_{B_{\mathcal{G}}}$ are irreducible.
\end{rem}
Using item (\emph{c}) in the Perron-Frobenius lemma, the spectral
radius $\rho\left(M\right)$ can be written as follows:
\begin{align}
\rho(M) & =\inf\{\lambda\colon\exists\mathbf{u}\succ\mathbf{0}\text{ s.t. }M\mathbf{u}\preceq\lambda\mathbf{u}\}\nonumber \\
 & =\inf\left\{ \lambda\colon\exists\mathbf{u}\succ\mathbf{0}\text{ s.t. }\max_{i\in\left[n\right]}\left\{ \sum_{j=1}^{n}M_{ij}\frac{u_{j}}{u_{i}}\right\} \leq\lambda\right\} \nonumber \\
 & =\inf\left\{ \lambda\colon\adjustlimits\inf_{\mathbf{u}\succ\mathbf{0}}\max_{i\in\left[n\right]}\left\{ \sum_{j=1}^{n}M_{ij}\frac{u_{j}}{u_{i}}\right\} \leq\lambda\right\} \nonumber \\
 & =\adjustlimits\inf_{\mathbf{u}\succ\mathbf{0}}\max_{i\in\left[n\right]}\left\{ \sum_{j=1}^{n}M_{ij}\frac{u_{j}}{u_{i}}\right\} ,\label{eq:Rho as inf-max}
\end{align}
where in the last equality we use the fact that $\inf\left\{ \lambda\colon a\leq\lambda\right\} =a$
for any $a$. Using (\ref{eq:Rho as inf-max}), we rewrite the optimization
problem~\eqref{eq:prob_robust_alloc} as
\begin{alignat}{2}
 & \underset{\mathbf{d}^{c}}{\optminc} & \quad & \adjustlimits\sup_{B_{\mathcal{G}}\in\Delta_{B_{\mathcal{G}}}}\inf_{\mathbf{u}\succ\mathbf{0}}\max_{i\in\left[n\right]}\left\{ \sum_{j=1}^{n}M_{ij}\left(B_{\mathcal{G}},\mathbf{d}^{c}\right)\frac{u_{j}}{u_{i}}\right\} \label{eq:prob_robust_alloc-2-1}\\
 & \optstc & \quad & \sum_{i\in\mathcal{V}_{C}}g_{i}\left(\delta_{i}^{c}\right)\leq C;\qquad\underline{\delta}_{i}^{c}\leq\delta_{i}^{c}\leq\overline{\delta}_{i}^{c},\quad\forall i\in\mathcal{V}_{C}.\nonumber 
\end{alignat}
In what follows, we will first cast problem~\eqref{eq:prob_robust_alloc-2-1}
into a robust GP with coefficient uncertainties. The main technical
challenge we face is the $\min$-$\sup$-$\inf$-$\max$ structure
in the objective function of problem \eqref{eq:prob_robust_alloc-2-1}.
As we prove in Appendix~\ref{sub:Proof-of-Proposition}, we can use
the Saddle Point Theorem (Proposition~\ref{prop:saddle_pt_thm},
Appendix \ref{sub:The-Saddle-Point-Theorem}) to exchange the order
of $\sup$ and $\inf$, so that problem~\eqref{eq:prob_robust_alloc-2-1}
can be written as
\begin{alignat*}{2}
 & \underset{\mathbf{d}^{c}}{\optminc} & \quad & \adjustlimits\inf_{\mathbf{u}\succ\mathbf{0}}\sup_{B_{\mathcal{G}}\in\Delta_{B_{\mathcal{G}}}}\max_{i\in\left[n\right]}\left\{ \sum_{j=1}^{n}M_{ij}\left(B_{\mathcal{G}},\mathbf{d}^{c}\right)\frac{u_{j}}{u_{i}}\right\} \\
 & \underset{\phantom{\mathbf{d}^{c}}}{\optstc} &  & \sum_{i\in\mathcal{V}_{C}}g_{i}\left(\delta_{i}^{c}\right)\leq C;\qquad\underline{\delta}_{i}^{c}\leq\delta_{i}^{c}\leq\overline{\delta}_{i}^{c},\quad\forall i\in\mathcal{V}_{C},
\end{alignat*}
which is equivalent to:
\begin{alignat*}{2}
 & \underset{\mathbf{d}^{c},\mathbf{u}\succ\mathbf{0}}{\optminc} & \quad & \adjustlimits\max_{i\in\left[n\right]}\sup_{B_{\mathcal{G}}\in\Delta_{B_{\mathcal{G}}}}\left\{ \sum_{j=1}^{n}M_{ij}\left(B_{\mathcal{G}},\mathbf{d}^{c}\right)\frac{u_{j}}{u_{i}}\right\} \\
 & \underset{\phantom{\mathbf{d}^{c},\mathbf{u}\succ\mathbf{0}}}{\optstc} & \quad & \sum_{i\in\mathcal{V}_{C}}g_{i}\left(\delta_{i}^{c}\right)\leq C;\qquad\underline{\delta}_{i}^{c}\leq\delta_{i}^{c}\leq\overline{\delta}_{i}^{c},\quad\forall i\in\mathcal{V}_{C}.
\end{alignat*}
If we introduce a slack variable 
\[
\lambda\triangleq\adjustlimits\max_{i\in\left[n\right]}\sup_{B_{\mathcal{G}}\in\Delta_{B_{\mathcal{G}}}}\left\{ \sum_{j=1}^{n}M_{ij}\left(B_{\mathcal{G}},\mathbf{d}^{c}\right)\frac{u_{j}}{u_{i}}\right\} ,
\]
we obtain the optimization problem described in the following proposition.
\begin{prop}
\label{prop:cgp_formulation}Assume $\mathcal{G}$ is a strongly connected
contact graph. The robust allocation problem~\eqref{eq:prob_robust_alloc}
achieves the same optimal value as the following optimization problem\emph{:}
\begin{alignat}{2}
 & \underset{\mathbf{d}^{c},\mathbf{u},\lambda}{\optmin}\quad &  & \lambda\label{eq:prob_robust_alloc_cgp}\\
 & \optst\quad &  & \sup_{B_{\mathcal{G}}\in\Delta_{B_{\mathcal{G}}}}\sum_{j=1}^{n}\beta_{ij}\frac{u_{j}}{u_{i}}+\delta_{i}^{c}\leq\lambda,\quad i\in\left[n\right],\nonumber \\
 &  &  & \sum_{i\in\mathcal{V}_{C}}g_{i}\left(\delta_{i}^{c}\right)\leq C;\qquad\underline{\delta}_{i}^{c}\leq\delta_{i}^{c}\leq\overline{\delta}_{i}^{c},\quad\forall i\in\mathcal{V}_{C},\nonumber \\
 &  &  & \mathbf{u}\succ\mathbf{0},\qquad\prod_{i=1}^{n}u_{i}=1.\label{eq:u_normalization}
\end{alignat}
Moreover, if $\mathbf{d}^{c*}$ optimizes problem~\eqref{eq:prob_robust_alloc_cgp},
then it is an optimal solution of problem~\eqref{eq:prob_robust_alloc}.\end{prop}
\begin{IEEEproof}
See Appendix~\ref{sub:Proof-of-Proposition}.
\end{IEEEproof}
Notice that the decision variables $\mathbf{d}^{c}$, $\mathbf{u}$,
and $\lambda$ of problem~\eqref{eq:prob_robust_alloc_cgp} are all
strictly positive. If $g_{i}$ is a monomial for all $i\in\mathcal{V}_{C}$,
then problem~\eqref{eq:prob_robust_alloc_cgp} is a robust GP with
coefficient uncertainties (in $\beta_{ij}$). In the next, we will
show that the uncertainty set $\Delta_{B_{\mathcal{G}}}$ can be relaxed
into a convex set in the form of~\eqref{eq:Tractable Set}, so that
problem~\eqref{eq:prob_robust_alloc_cgp} can be solved efficiently
as a conic GP using the methodology proposed in Section~\ref{sub:Robust-Geometric-Programming}.

%% file: uncertainty_set.tex
As we mentioned in Section~\ref{sub:Problem-Formulation}, the uncertainty
set $\Delta_{B_{\mathcal{G}}}^{D}$ is nonconvex since it is defined
by a collection of polynomial equalities. In this subsection, we define
a convex superset $\widehat{\Delta}_{B_{\mathcal{G}}}^{D}\supset\Delta_{B_{\mathcal{G}}}^{D}$,
so that problem~(\ref{eq:prob_robust_alloc}) becomes a conic geometric
program after we substitute $\Delta_{B_{\mathcal{G}}}^{D}$ by $\widehat{\Delta}_{B_{\mathcal{G}}}^{D}$
(which changes the combined uncertainty set $\Delta_{B_{\mathcal{G}}}$). 

We define the convex superset $\widehat{\Delta}_{B_{\mathcal{G}}}^{D}$
as follows: 
\begin{multline}
\widehat{\Delta}_{B_{\mathcal{G}}}^{D}\triangleq\Biggl\{ B_{\mathcal{G}}\in\mathbb{R}^{n\times n}\colon\frac{1}{n}\sum_{j\in\mathcal{V}_{S}}\beta_{ij}p_{j}(t)\\
\leq1-\left(1-\frac{p_{i}(t+1)-p_{i}(t)\left(1-\delta_{i}^{0}\right)}{1-p_{i}(t)}\right)^{1/n}\\
\mbox{for all }i\in\mathcal{V}_{S},\ t\in\left[T\right]\text{ s.t. }p_{i}(t)<1\Biggr\}.\label{eq:Superset}
\end{multline}

\begin{lem}
\label{lem:Superset}The set~$\widehat{\Delta}_{B_{\mathcal{G}}}^{D}$
is a superset of~$\Delta_{B_{\mathcal{G}}}^{D}$.\end{lem}
\begin{IEEEproof}
Consider any $i\in\mathcal{V}_{S}$. Recall that it always holds that
$1-\beta_{ij}p_{j}(t)\geq0$. Then, from the AM-GM inequality, we
have that
\begin{align}
\prod_{j=1}^{n}[1-\beta_{ij}p_{j}(t)] & \leq\left(\frac{\sum_{j=1}^{n}[1-\beta_{ij}p_{j}(t)]}{n}\right)^{n}\nonumber \\
 & =\left(1-\frac{1}{n}\sum_{j=1}^{n}\beta_{ij}p_{j}(t)\right)^{n}.\label{eq:prod_ineq}
\end{align}
We can use~(\ref{eq:prod_ineq}) to yield a constraint on the transmission
rates $\beta_{ij}$ from the empirical dataset~$\mathcal{D}$ defined
in~(\ref{eq:Data Series}). Applying~(\ref{eq:prod_ineq}) to the
nonlinear dynamics~(\ref{eq:HeNiSIS dynamics}) results in the following
inequality: 
\begin{multline*}
p_{i}(t+1)\geq p_{i}(t)(1-\delta_{i}^{0})\\
+(1-p_{i}(t))\left\{ 1-\left(1-\frac{1}{n}\sum_{j=1}^{n}\beta_{ij}p_{j}(t)\right)^{n}\right\} ,
\end{multline*}
where we have used the fact that the recovery rate $\delta_{i}$ is
equal to the natural recovery rate $\delta_{i}^{0}$ during empirical
observations. Since $p_{i}(t)<1$, we can rearrange the above inequality
to obtain
\[
\left(1-\frac{1}{n}\sum_{j=1}^{n}\beta_{ij}p_{j}(t)\right)^{n}\geq1-\frac{p_{i}(t+1)-p_{i}(t)(1-\delta_{i}^{0})}{1-p_{i}(t)},
\]
which is equivalent to
\begin{equation}
1-\frac{1}{n}\sum_{j=1}^{n}\beta_{ij}p_{j}(t)\geq\left(1-\frac{p_{i}(t+1)-p_{i}(t)(1-\delta_{i}^{0})}{1-p_{i}(t)}\right)^{1/n}.\label{eq:am-gm_2}
\end{equation}
Here we have used the fact 
\[
1-\frac{p_{i}(t+1)-p_{i}(t)(1-\delta_{i}^{0})}{1-p_{i}(t)}=\prod_{j=1}^{n}[1-\beta_{ij}p_{j}(t)]>0
\]
according to the system dynamics~(\ref{eq:HeNiSIS dynamics}). We
can rearrange~(\ref{eq:am-gm_2}) to obtain
\[
\frac{1}{n}\sum_{j=1}^{n}\beta_{ij}p_{j}(t)\leq1-\left(1-\frac{p_{i}(t+1)-p_{i}(t)(1-\delta_{i}^{0})}{1-p_{i}(t)}\right)^{1/n}.
\]
Finally, we use the fact 
\[
\frac{1}{n}\sum_{j\in\mathcal{V}_{S}}\beta_{ij}p_{j}(t)\leq\frac{1}{n}\sum_{j=1}^{n}\beta_{ij}p_{j}(t)
\]
to complete the proof.
\end{IEEEproof}
The following comments are in order. First, the uncertainty set $\widehat{\Delta}_{B_{\mathcal{G}}}^{D}$
in (\ref{eq:Superset}) is defined by a collection of affine inequalities;
therefore, it is a convex polytope and can be represented in the form
of (\ref{eq:Tractable Set}). Second, from Lemma \ref{lem:Superset},
we have that the superset $\widehat{\Delta}_{B_{\mathcal{G}}}^{D}$
contains all the contact matrices $B_{\mathcal{G}}$ that are coherent
with both prior information and empirical observations. In the following
section, we illustrate our relaxation with numerical simulations and
verify that the robust allocation is not overly conservative. In fact,
for some realistic cases, the robust allocation achieves similar performance
as the optimal allocation solved under a known contact matrix $B_{\mathcal{G}}$.

\textbf{}

%% file: simulations.tex
\textcolor{red}{}In this section, we illustrate the robust data-driven
allocation framework developed in Section \ref{sec:Data-Driven-Resource-Allocation}.
We consider the problem of controlling an epidemic outbreak propagating
through the worldwide air transportation network~\cite{Schneider11}.
The nodes in the network represent airports, whereas edges are flight
connections for which we know passenger flows. Through our simulations,
we demonstrate the following facts about the data-driven allocation
framework. First, incorporating observations into the uncertainty
set $\Delta_{B\mathcal{_{G}}}$ significantly helps reduce the worst-case
spreading rate bound. Second, the robust allocation framework does
not need many observations to converge; in particular, the length
of the observation period that we need is only a fraction of the number
of nodes in the network. Finally, even though the robust allocation
algorithm does not have access to the true underlying contact network~$B_{\mathcal{G}}$,
the resulting allocation achieves very similar performance compared
to the optimal allocation solved using the actual $B_{\mathcal{G}}$.

\subsection{Numerical Setup}

\begin{figure}
\begin{centering}
\includegraphics[width=2.81in]{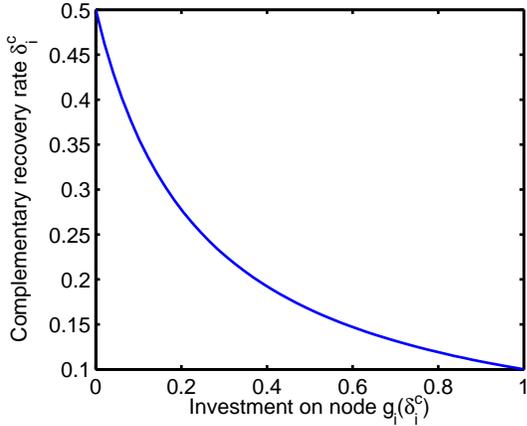}
\par\end{centering}

\caption{Plot of the inverse of the vaccination cost function $g_{i}^{-1}$.
This function represents the complementary recovery rate $\delta_{i}^{c}$
as a function of the investment on that node.}
\label{fig:VaccineCost}
\end{figure}

In our simulations, we consider the problem of controlling an epidemic
outbreak propagating through a flight network comprised by the top
$100$ airports (based on yearly total traffic), so that $n=100$.
To illustrate the robust data-driven approach, we first generate a
time series representing the dynamics of a hypothetical outbreak using
the nonlinear dynamics~(\ref{eq:HeNiSIS dynamics}). We run our simulation
assuming a homogeneous value for the natural recovery rate, $\overline{\delta}_{i}^{c}=0.5$
for all $i\in\left[n\right]$, and a link-dependent contact rate $\beta_{ij}$
that is proportional to the traffic through that edge. Assuming an
initial infection $p_{i}(0)=0.5$ for all $i\in[n]$, we generate
a time series $\{\mathbf{p}(t)\}_{t=1}^{T}$ representing the evolution
of the infection over time.

In our data-driven framework, we assume that we do not have direct
access to the matrix of infection rates $B_{\mathcal{G}}$. Instead,
the data-driven allocation algorithm only has access to the observations
$\{\mathbf{p}(t)\}_{t=1}^{T}$ for some period $t\in\left[T\right]$.
Using this data, our algorithm generates an uncertainty set $\Delta_{B_{\mathcal{G}}}$
of data-coherent contact matrices. The parameters that define the
uncertainty set $\Delta_{B_{\mathcal{G}}}$ are chosen as follows.
For all $\left(i,j\right)\in\mathcal{E}$, we assume an \emph{a prior}i
upper bound $\overline{\beta}_{ij}=1.5\beta_{ij}$ (i.e., the contact
rate of an edge is at most 50\% above its nominal contact rate), whereas
the lower bound is $\underline{\beta}_{ij}=0.5\beta_{ij}$ (i.e.,
the contact rate is at least half the nominal value). The natural
recovery rate is chosen as $\overline{\delta}_{i}^{c}=0.5$ for all
$i\in[n]$, while the lower bound is chose to be $\underline{\delta}_{i}^{c}=0.1$
for all $i\in[n]$ (i.e., the recovery rate $\delta_{i}$ is at most
$1-\underline{\delta}_{i}^{c}=0.9$). We assume that the set of control
nodes $\mathcal{V}_{C}=[n]$ and vary the set of sensing nodes.

\begin{figure}
\begin{centering}
\includegraphics[width=2.87in]{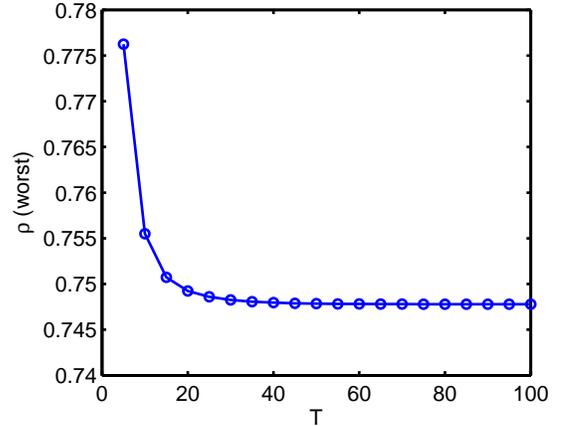}
\par\end{centering}

\caption{Evolution of the worst-case spectral radius $\rho_{\mathrm{wor}}$
as a function of the number of observations~$T$.}
\label{fig:rho_worst_T}
\end{figure}

To find the optimal allocation of vaccines, we consider the following
vaccination cost function $g_{i}$ given by 
\[
g_{i}(\delta_{i}^{c})=\frac{1/\delta_{i}^{c}-1/\overline{\delta}_{i}^{c}}{1/\underline{\delta}_{i}^{c}-1/\overline{\delta}_{i}^{c}}
\]
for all~$i\in[n]$. It can be seen that $g_{i}$ is a monomial in
$\delta_{i}^{c}$. The function $g_{i}$ satisfies $g_{i}(\overline{\delta}_{i}^{c})=0$;
namely, there is no cost by keeping $\delta_{i}^{c}$ as the natural
complementary recovery rate~$\overline{\delta}_{i}^{c}$. This function
also satisfies $g_{i}(\underline{\delta}_{i}^{c})=1$; namely, the
maximum allocation per node is one unit. Furthermore, the cost function
is monotonically decreasing and exhibits diminishing returns (see
Fig. \ref{fig:VaccineCost}). In this setup, our problem is to find
the optimal allocation of vaccines throughout the airports in the
air transportation network assuming we have a total budget equal to
$C=0.5n=50$.

\subsection{Results and Discussions}

For any given uncertainty set $\Delta_{B_{\mathcal{G}}}$, we define
the worst-case spectral radius~$\rho_{\mathrm{wor}}(\Delta_{B_{\mathcal{G}}})$
as the optimal value of the robust allocation problem~(\ref{eq:prob_robust_alloc}).
In other words, $\rho_{\mathrm{wor}}(\Delta_{B_{\mathcal{G}}})$ represents
the slowest exponential rate of disease eradication that can be achieved
for those contact matrices that are coherent with our observations.
In our first experiment, we illustrate the dependency of $\rho_{\mathrm{wor}}$
with respect to~$T$, i.e., the number of observations available.
In Fig.~\ref{fig:rho_worst_T}, we show the value of $\rho_{\mathrm{wor}}$
as we increase the observation period in the range $T=1,\ldots,100$.
Notice how, as $T$ grows, the amount of available information about
the contact network increases and, as a result, $\rho_{\mathrm{wor}}$
decreases (i.e., we are able to guarantee a faster disease eradication).
Notice also that the value of $\rho_{\mathrm{wor}}$ remains approximately
unchanged after $T=30$ observations. This result may seem surprising
at first glance, since from the perspective of system observability,
one would normally need as many time steps as the dimension of the
system (in this case, $n=100$) in order to identify the system. This
demonstrates one of the benefits of using the robust allocation framework;
namely, \emph{it allows us to find an allocation without performing
a previous system identification}.

\begin{figure}
\begin{centering}
\includegraphics[width=2.87in]{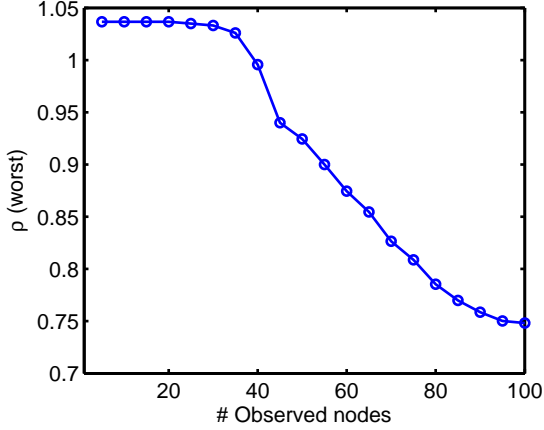}
\par\end{centering}

\caption{Evolution of the worst-case spectral radius $\rho_{\mathrm{wor}}$
as a function of the number of observed nodes $\left|\mathcal{V}_{S}\right|$.}
\label{fig:rho_worst_n_obs}
\end{figure}

In a second set of experiments, we numerically verify the performance
of our data-driven allocation algorithm in the presence of sparse
observations. In particular, we assume that we can only measure the
evolution of the disease in a set of sensor nodes $\mathcal{V}_{S}$,
which we choose to be those airports with the highest yearly total
traffic. In Fig.~\ref{fig:rho_worst_n_obs}, we plot the value of
$\rho_{\mathrm{wor}}$ as we increase the number of sensor nodes from
$\left|\mathcal{V}_{S}\right|=1,\ldots,100$. Notice how, as we increase
the number of sensor nodes, $\rho_{\mathrm{wor}}$ decreases. Interestingly,
for $\left|\mathcal{V}_{S}\right|\leq20$ sensors, the value of $\rho_{\mathrm{wor}}$
hardly changes. In contrast, we observe a dramatic improvement in
the value of $\rho_{\mathrm{wor}}$ for $\left|\mathcal{V}_{S}\right|\geq40$.
In fact, using only 40 sensors (out of 100 nodes), we can find an
allocation that guarantees the eradication of the disease (i.e., $\rho_{\mathrm{wor}}<1$),
even for the worst instantiation of $B_{\mathcal{G}}$ in $\Delta_{B_{\mathcal{G}}}$.

\begin{figure}
\begin{centering}
\includegraphics[width=2.87in]{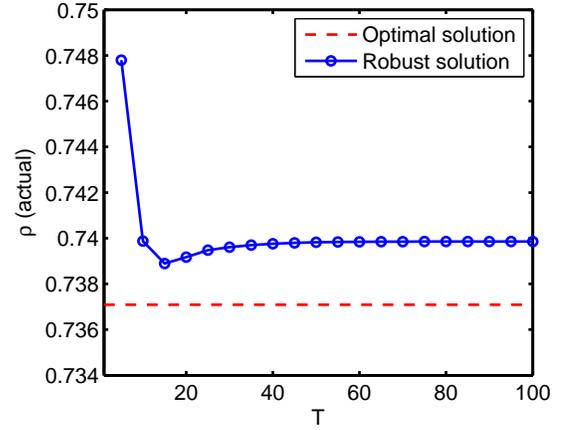}
\par\end{centering}

\caption{The (actual) spectral radius $\rho(M(B_{\mathcal{G}},\mathbf{d}^{c}))$
corresponding to both the optimal allocation $\mathbf{d}_{\mathrm{opt}}^{c}$
and the robust allocation $\mathbf{d}_{\mathrm{rob}}^{c}\left(T\right)$.
Note that it is not guaranteed that $\rho(M(B_{\mathcal{G}},\mathbf{d}_{\mathrm{rob}}^{c}))$
decreases monotonically with $T$. }
\label{fig:rho_actual}
\end{figure}

In our final simulation, we compare the allocation obtained from the
data-driven framework with the allocation obtained assuming we have
full access to the actual matrix of infection rates $B_{\mathcal{G}}$.
Using the framework proposed in Preciado et al.~\cite{PZEJP14},
we can obtain the optimal allocation $\mathbf{d}_{\mathrm{opt}}^{c}$,
which is defined as the solution to the following optimization problem:
\begin{alignat}{2}
 & \underset{\mathbf{\mathbf{d}^{c}}}{\optmin} & \quad & \rho(M(B_{\mathcal{G}},\mathbf{d}^{c}))\label{eq:prob_opt_alloc}\\
 & \optst & \quad & \sum_{i\in\mathcal{V}_{C}}g_{i}\left(\delta_{i}^{c}\right)\leq C,\nonumber \\
 &  &  & \delta_{i}^{c}\leq\delta_{i}^{c}\leq\overline{\delta}_{i}^{c},\quad i\in\mathcal{V}_{C}.\nonumber 
\end{alignat}
The optimal value $\rho(M(B_{\mathcal{G}},\mathbf{d}_{\mathrm{opt}}^{c}))$
of problem~(\ref{eq:prob_opt_alloc}) represents the fastest exponential
rate at which the disease is eradicated when the contact network is
completely known. Additionally, we denote by $\mathbf{d}_{\mathrm{rob}}^{c}\left(T\right)$
the optimal solution to the robust data-driven allocation problem~(\ref{eq:prob_robust_alloc})
when $T$ time samples are available. We evaluate the spectral radius~$\rho(M(B_{\mathcal{G}},\mathbf{d}_{\mathrm{rob}}^{c}\left(T\right)))$,
which represents the exponential rate at which the disease is eradicated
when we apply the allocation $\mathbf{d}_{\mathrm{rob}}^{c}\left(T\right)$
to the actual contact network $B_{\mathcal{G}}$. In Fig.~\ref{fig:rho_actual},
we compare $\rho(M(B_{\mathcal{G}},\mathbf{d}_{\mathrm{rob}}^{c}\left(T\right)))$
with the optimal value $\rho(M(B_{\mathcal{G}},\mathbf{d}_{\mathrm{opt}}^{c}))$
for different values of $T$. Since $\mathbf{d}_{\mathrm{opt}}^{c}$
is the optimal solution to problem~(\ref{eq:prob_opt_alloc}), we
always have $\rho(M(B_{\mathcal{G}},\mathbf{d}_{\mathrm{opt}}^{c}))\leq\rho(M(B_{\mathcal{G}},\mathbf{d}_{\mathrm{rob}}^{c}))$.
However, Fig.~\ref{fig:rho_actual} shows that the difference between
$\rho(M(B_{\mathcal{G}},\mathbf{d}_{\mathrm{opt}}^{c}))$ and $\rho(M(B_{\mathcal{G}},\mathbf{d}_{\mathrm{rob}}^{c}\left(T\right)))$
is small for the particular network under investigation.

Finally, it is worth mentioning that the robust data-driven allocation
problem does not take significantly more time to solve than the optimal
allocation problem. We have solved both allocation problems in MATLAB
(R2012b) using CVX (Version 2.1, Build 1079)~\cite{cvx} with the
Mosek solver (Version 7.0.0.106). All computations are carried out
on a laptop computer equipped with a dual-core 2.5 GHz Intel Core
i5 processor and 4 GB of RAM. For $n=100$, the optimal allocation
problem takes approximately 17 seconds to solve, whereas the robust
allocation problem with \emph{a priori} bounds on $\beta_{ij}$ takes
approximately 49 seconds for $T=30$.

%% file: conclusions.tex
We have introduced a novel mathematical framework, based on conic
geometric programming, to control a viral spreading process taking
place in a contact network with unknown contact rates. We assume that
we have access to time series data describing the evolution of the
spreading process over a finite time period over a collection of sensor
nodes. Using this data, we have developed a data-driven robust convex
optimization framework to find the optimal allocation of protection
resources over a set of control nodes to eradicate the viral spread
at the fastest possible rate.

We have illustrated our approach using data obtained from the worldwide
air transportation network. We have simulated a hypothetical epidemic
outbreak over a finite time period and fed the resulting time series
in our data-driven optimization algorithm. From our numerical results,
we verify that (\emph{i}) incorporating observations into the data-driven
allocation algorithm significantly reduces the worst-case spreading
rate bound; (\emph{ii}) the robust allocation framework does not need
many observations to converge; (\emph{iii}) even though the robust
allocation algorithm does not have access to the true underlying contact
network~$B_{\mathcal{G}}$, the resulting allocation achieves very
similar performance compared to the optimal allocation solved under
the actual $B_{\mathcal{G}}$.

%% file: saddle_pt_thm.tex
The purpose of this section is to present the Saddle Point Theorem
from convex analysis that is used in the proof of Proposition~\ref{prop:cgp_formulation}.
Most of this section is adopted from Chapters 1 and 2 of the book
by Bertsekas et al.~\cite{bertsekas2003convex}. To prepare for this,
we first introduce some basic concepts in convex analysis. 
\begin{defn}
[Epigraph]Let $X$ be a subset of $\mathbb{R}^{n}$. The \emph{epigraph}
of an extended real-valued function $f\colon X\to[-\infty,\infty]$
is defined as the set 
\[
\mathrm{epi}(f)=\{(x,w\}\colon x\in X,\; w\in\mathbb{R},\; f(x)\leq w\}.
\]

\end{defn}
\begin{defn}
[Closed Function]Let $X$ be a subset of $\mathbb{R}^{n}$. An extend
real-valued function $f\colon X\to[-\infty,\infty]$ is called \emph{closed}
if its epigraph $\mathrm{epi}(f)$ is a closed set. 
\end{defn}
\begin{defn}
[Convex Function]Let $C$ be a convex subset of $\mathbb{R}^{n}$.
An extended real-valued function $f\colon C\to[-\infty,\infty]$ is
called \emph{convex} if $\mathrm{epi}(f)$ is a convex set. 
\end{defn}
Let $X$ and $Z$ be nonempty convex subsets of $\mathbb{R}^{n}$
and $\mathbb{R}^{m}$, respectively. The Saddle Point Theorem considers
a real-valued function $\phi\colon X\times Z\to\mathbb{R}$ and provides
conditions to ascertain the \emph{minimax equality}
\begin{equation}
\adjustlimits\sup_{z\in Z}\inf_{x\in X}\phi(x,z)=\adjustlimits\inf_{x\in X}\sup_{z\in Z}\phi(x,z).\label{eq:minimax_eq}
\end{equation}
Before introducing the Saddle Point Theorem, for each $z\in Z$, we
define the function $t_{z}\colon\mathbb{R}^{n}\to(-\infty,\infty]$
as
\[
t_{z}(x)=\begin{cases}
\phi(x,z) & \quad x\in X\\
\infty & \quad x\notin X,
\end{cases}
\]
and, for each $x\in X$, we define the function $r_{x}\colon\mathbb{R}^{m}\to(-\infty,\infty]$
as
\[
r_{x}(z)=\begin{cases}
-\phi(x,z) & \quad z\in Z\\
\infty & \quad z\notin Z.
\end{cases}
\]
We also need the following assumption on $t_{z}$ and $r_{x}$ (or
equivalently, $\phi$).

\begin{assumption}

\label{asu:closed_and_convex}The function~$t_{z}$ is closed and
convex for each $z\in Z$, and the function $r_{x}$ is closed and
convex for each $x\in X$.

\end{assumption}
\begin{rem}
\label{rem:sufficient_conditions_rt}One useful sufficient condition
for the function $t_{z}$ to be closed is that the set $X$ is closed
and the function $\phi(x,z)$ is lower semicontinous in~$x$. In
addition, the convexity of $t_{z}$ is equivalent to the convexity
of $\phi(x,z)$ in $x$ over $X$. A similar sufficient condition
can also be applied to $r_{x}$. 
\end{rem}
We are now ready to present the Saddle Point Theorem.
\begin{prop}
[Saddle Point Theorem]\label{prop:saddle_pt_thm}Suppose $\phi(x,z)$
satisfies Assumption~\ref{asu:closed_and_convex}. Then $\phi$ satisfies
the minimax equality~(\ref{eq:minimax_eq}) under any of the following
conditions\end{prop}
\begin{enumerate}
\item $X$ and $Z$ are compact.
\item $Z$ is compact, and there exists $\bar{z}\in Z$ and $\gamma\in\mathbb{R}$
such that the set 
\[
\{x\in X\colon\phi(x,\bar{z})\leq\gamma\}
\]
is nonempty and compact.
\item $X$ is compact, and there exists $\bar{x}\in X$ and $\gamma\in\mathbb{R}$
such that the set 
\[
\{z\in Z\colon\phi(\bar{x},z)\geq\gamma\}
\]
is nonempty and compact.
\item There exist $\bar{x}\in X$, $\bar{z}\in Z$, and $\gamma\in\mathbb{R}$
such that the sets
\[
\{x\in X\colon\phi(x,\bar{z})\leq\gamma\},\qquad\{z\in Z\colon\phi(\bar{x},z)\geq\gamma\}
\]
are nonempty and compact.\end{enumerate}

%% file: proof_cgp_formulation.tex
Without loss of generality, we assume that the vector $\mathbf{u}$
satisfies $\prod_{i=1}^{n}u_{i}=1$. Define the set
\[
\mathcal{U}\triangleq\left\{ \mathbf{u}\in\mathbb{R}^{n}\colon\mathbf{u}\succ\mathbf{0},\ \prod_{i=1}u_{i}=1\right\} .
\]
We wish to show that the following minimax equality holds:
\begin{multline}
\adjustlimits\sup_{B_{\mathcal{G}}\in\Delta_{B_{\mathcal{G}}}}\inf_{\mathbf{u}\in\mathcal{U}}\max_{i\in\left[n\right]}\left\{ \sum_{j=1}^{n}M_{ij}\left(B_{\mathcal{G}},\mathbf{d}^{c}\right)u_{j}/u_{i}\right\} \\
=\adjustlimits\inf_{\mathbf{u}\in\mathcal{U}}\sup_{B_{\mathcal{G}}\in\Delta_{B_{\mathcal{G}}}}\max_{i\in\left[n\right]}\left\{ \sum_{j=1}^{n}M_{ij}\left(B_{\mathcal{G}},\mathbf{d}^{c}\right)u_{j}/u_{i}\right\} .\label{eq:minimax_lambda1}
\end{multline}
Define $\tilde{u}_{i}=\log u_{i}$ for all $i\in[n]$ and the function
\begin{align*}
\phi(\widetilde{\mathbf{u}},B_{\mathcal{G}}) & =\max_{i\in\left[n\right]}\left\{ \sum_{j=1}^{n}M_{ij}\left(B_{\mathcal{G}},\mathbf{d}^{c}\right)\exp(\widetilde{u}_{j}-\widetilde{u}_{i})\right\} \\
 & =\max_{i\in\left[n\right]}\left\{ \sum_{j=1}^{n}\beta_{ij}\exp(\widetilde{u}_{j}-\widetilde{u}_{i})+\delta_{i}^{c}\right\} .
\end{align*}
Then, the minimax equality~(\ref{eq:minimax_lambda1}) holds if and
only if the following equality holds:
\begin{equation}
\adjustlimits\sup_{B_{\mathcal{G}}\in\Delta_{B_{\mathcal{G}}}}\inf_{\widetilde{\mathbf{u}}\in\widetilde{\mathcal{U}}}\phi(\widetilde{\mathbf{u}},B_{\mathcal{G}})=\adjustlimits\inf_{\widetilde{\mathbf{u}}\in\widetilde{\mathcal{U}}}\sup_{B_{\mathcal{G}}\in\Delta_{B_{\mathcal{G}}}}\phi(\widetilde{\mathbf{u}},B_{\mathcal{G}}),\label{eq:minimax_lambda1_exp}
\end{equation}
where $\mathcal{\widetilde{U}}\triangleq\left\{ \widetilde{\mathbf{u}}\in\mathbb{R}^{n}\colon\sum_{i=1}^{n}\widetilde{u}_{i}=0\right\} $.

Using the sufficient conditions in Remark~\ref{rem:sufficient_conditions_rt},
it can be verified that $\phi$ satisfies Assumption~\ref{asu:closed_and_convex}.
To apply the Saddle Point Theorem (Proposition~\ref{prop:saddle_pt_thm}),
we substitute $x=\widetilde{\mathbf{u}}$, $X=\widetilde{\mathcal{U}}$,
$z=B_{\mathcal{G}}$, and $Z=\Delta_{B_{\mathcal{G}}}$ in the Saddle
Point Theorem. Since $\Delta_{B_{\mathcal{G}}}$ is compact according
to item (\emph{i}) in the formulation of Problem~\ref{prob:Main Problem},
we can apply the second condition in the Saddle Point Theorem if we
can show that there exist $\widehat{B}_{\mathcal{G}}\in\Delta_{B_{\mathcal{G}}}$
and $\gamma\in\mathbb{R}$ such that the set
\begin{equation}
\mathcal{S}_{\widetilde{\mathbf{u}}}\triangleq\left\{ \widetilde{\mathbf{u}}\in\widetilde{\mathcal{U}}\colon\phi(\widetilde{\mathbf{u}},\widehat{B}_{\mathcal{G}})\leq\gamma\right\} \label{eq:def_Su}
\end{equation}
is nonempty and compact.

Consider any $\widehat{B}_{\mathcal{G}}\in\Delta_{B_{\mathcal{G}}}$
and choose $\gamma=\phi(\mathbf{0},\widehat{B}_{\mathcal{G}})$. We
can observe that~$\mathcal{S}_{\widetilde{\mathbf{u}}}$ is nonempty,
since $\mathbf{0}\in\widetilde{\mathcal{U}}$. In order to show that
$\mathcal{S}_{\widetilde{\mathbf{u}}}$ is compact, notice that $\mathcal{S}_{\widetilde{\mathbf{u}}}$
is a subset of $\mathbb{R}^{n}$. Then, compactness of $\mathcal{S}_{\widetilde{\mathbf{u}}}$
is equivalent to that $\mathcal{S}_{\widetilde{\mathbf{u}}}$ is closed
and bounded as a result of the Heine\textendash{}Borel theorem. To
show that $\mathcal{S}_{\widetilde{\mathbf{u}}}$ is closed, we use
the fact that $\phi$ is continuous in $\widetilde{\mathbf{u}}$,
and $\mathcal{S}_{\widetilde{\mathbf{u}}}$ is the preimage of the
closed set $\{y\in\mathbb{R}\colon y\leq\gamma\}$. To show that $\mathcal{S}_{\widetilde{\mathbf{u}}}$
is bounded, suppose by contradiction that $\mathcal{S}_{\widetilde{\mathbf{u}}}$
is unbounded, which implies that there exists $\widehat{\mathbf{u}}\in\mathcal{S}_{\widetilde{\mathbf{u}}}$
such that $\widehat{\mathbf{u}}\neq\mathbf{0}$ and $\alpha\widehat{\mathbf{u}}\in\mathcal{S}_{\widetilde{\mathbf{u}}}$
for all $\alpha>0$. Substituting $\alpha\widehat{\mathbf{u}}\in\mathcal{S}_{\widetilde{\mathbf{u}}}$
into~(\ref{eq:def_Su}), we obtain
\[
\phi(\alpha\widehat{\mathbf{u}},\widehat{B}_{\mathcal{G}})=\max_{i\in\left[n\right]}\left\{ \sum_{j=1}^{n}\widehat{\beta}_{ij}\exp(\alpha(\widehat{u}_{j}-\widehat{u}_{i}))+\delta_{i}^{c}\right\} \leq\gamma
\]
for all $\alpha>0$. Define $\mathcal{I}^{*}\triangleq\{i\in[n]\colon\widehat{u}_{i}=\max_{i\in[n]}\{\widehat{u}_{i}\}\}$.
We know that the set $[n]\backslash\mathcal{I}^{*}$ is nonempty;
otherwise we have $\widehat{u}_{1}=\widehat{u}_{2}=\dots=\widehat{u}_{n}$
and hence $\widehat{\mathbf{u}}=\mathbf{0}$. Then, from the irreducibility
of~$\widehat{B}_{\mathcal{G}}$, we know that there exist $k\in\mathcal{I}^{*}$
and $\ell\in[n]\backslash\mathcal{I}^{*}$ such that $\widehat{\beta}_{k\ell}\neq0$.
Using the definition of $\mathcal{I}^{*}$, we know that $\widehat{u}_{k}-\widehat{u}_{\ell}>0$.
As $\alpha\to\infty$, we have $\widehat{\beta}_{k\ell}\exp(\alpha(\widehat{u}_{k}-\widehat{u}_{\ell}))\to\infty$
and hence $\phi(\alpha\widehat{\mathbf{u}},\widehat{B}_{\mathcal{G}})\to\infty$,
which leads to a contradiction. 

To summarize, we have shown that~$\Delta_{B_{\mathcal{G}}}$ is compact
and there exists $\gamma\in\mathbb{R}$ such that $\mathcal{S}_{\widetilde{\mathbf{u}}}$
is nonempty and compact. Then, from the second condition in the Saddle
Point Theorem, we know that the minimax equality~(\ref{eq:minimax_lambda1_exp})
holds. This allows us to rewrite the objective of problem~(\ref{eq:prob_robust_alloc-2-1})
as
\[
\adjustlimits\inf_{\mathbf{u}\in\mathcal{U}}\sup_{B_{\mathcal{G}}\in\Delta_{B_{\mathcal{G}}}}\max_{i\in\left[n\right]}\left\{ \sum_{j=1}^{n}M_{ij}\left(B_{\mathcal{G}},\mathbf{d}^{c}\right)u_{j}/u_{i}\right\} .
\]
By introducing a slack variable 
\[
\lambda\triangleq\adjustlimits\max_{i\in\left[n\right]}\sup_{B_{\mathcal{G}}\in\Delta_{B_{\mathcal{G}}}}\left\{ \sum_{j=1}^{n}M_{ij}\left(B_{\mathcal{G}},\mathbf{d}^{c}\right)u_{j}/u_{i}\right\} ,
\]
we obtain the optimization problem described in Proposition~\ref{prop:cgp_formulation},
where the constraint~(\ref{eq:u_normalization}) is given by the
definition of~$\mathcal{U}$.